\documentclass[11pt]{amsart}

\usepackage{pb-diagram}
\usepackage{amssymb, stmaryrd}

\theoremstyle{plain}
\newtheorem{theorem}{Theorem}[section]

\newtheorem{lemma}[theorem]{Lemma}

\numberwithin{equation}{section}

\theoremstyle{definition}

\theoremstyle{remark}
\newtheorem{remark}[theorem]{Remark}

\newcommand{\la}{\langle}
\newcommand{\ra}{\rangle}

\newcommand{\ow}{\omega}

\newcommand{\R}{\mathbb{R}}
\newcommand{\C}{\mathbb{C}}

\newcommand{\fg}{\mathfrak{g}}

\newcommand{\U}{\operatorname{U}}
\newcommand{\Su}{\operatorname{SU}}
\newcommand{\Sp}{\operatorname{Sp}}
\renewcommand{\H} {\operatorname{H}}

\DeclareMathOperator {\PP} {\mathbb{P}}
\DeclareMathOperator {\CP} {\mathbb{CP}}

\newcommand     {\comment}[1]   {}
\newcommand{\mute}[2] {}
\newcommand     {\printname}[1] {}

\newcommand{\labell}[1] {\label{#1}\printname{#1}}

\begin{document}

\title{Strong Lefschetz property under reduction}
\author{River Chiang}
\address{Department of Mathematics, National Cheng Kung University,
Tainan 701, Taiwan}
\email{riverch@mail.ncku.edu.tw}
\author{Eugene Z.~Xia}
\email{ezxia@mail.ncku.edu.tw}
\begin{abstract}
Let $n>1$ and \(G\) be the group \(\Su(n)\) or \(\Sp(n)\). This
paper constructs compact symplectic manifolds whose symplectic
quotient under a Hamiltonian \(G\)-action does not inherit the
strong Lefschetz property.
\end{abstract}
\maketitle
\date{today}

\section{Introduction}\labell{S:intro}

Let $(M, \Omega)$ be a compact symplectic manifold of dimension $2m$.
For each $0 \leq k \leq m$, the symplectic form $\Omega$ induces a map
\[
\Omega^{k} \colon \H^{m-k}(M) \to \H^{m+k}(M), \qquad [\alpha]
\mapsto [\Omega^{k}\wedge\alpha].
\]
The manifold \((M,\Omega)\) is called Lefschetz if \(\Omega^{m-1}\)
is an isomorphism.  It is called strong Lefschetz
if \(\Omega^k\) is an isomorphism for all \(0 \le k \le m\).
By Poincar\'e Duality, the Lefschetz maps are surjective if and
only if they are injective.

The simplest examples of symplectic manifolds are K\"ahler
manifolds. By the Strong Lefschetz Theorem, all compact K\"ahler
manifolds are strong Lefschetz. On the other hand, Thurston and
others constructed compact symplectic manifolds that do not admit
K\"ahler structures \cite{bg,cfg,gompf,mcduff,thurston}, either by
showing that these manifolds have an odd Betti number which is
odd, or that they do not have the strong Lefschetz property.

Another interesting fact about the Lefschetz property concerns
Hamiltonian actions. A Theorem of Ono \cite{ono} asserts that
if a symplectic manifold $(M, \Omega)$ is
Lefschetz, then a symplectic circle action on $(M,
\Omega)$ is Hamiltonian if and only if it has fixed points.

We recall that an action of a compact Lie group $G$ on $(M,
\Omega)$ is symplectic if it preserves the symplectic form $\Omega$.
Let $\fg$ be the Lie algebra of $G$ and $\fg^*$ its dual. For any
$\xi\in \fg$, denote by $\xi_M$ the vector field
on $M$ induced by the group action of $\exp (t\xi)$.
A symplectic action is Hamiltonian if there exists a moment map
$\Phi\colon M \to \fg^*$,
which is equivariant with respect to the group action on $M$
and the coadjoint action on $\fg^*$ such that
\[
d\la \Phi, \xi\ra = \iota(\xi_M)\Omega
\] for all $\xi\in \fg$. The orbit space $\Phi^{-1}(0)/G$ is called the
symplectic quotient, and often denoted by $M{\sslash}G$.
If $G$ acts freely on $\Phi^{-1}(0)$, the symplectic quotient $M{\sslash}G$ is
a smooth symplectic manifold. Its symplectic form $\ow$ is defined by
\[
\pi^*\ow=i^*\Omega,
\]
where $\pi\colon \Phi^{-1}(0) \to M{\sslash}G$ is the orbit map,
and $i\colon \Phi^{-1}(0)\hookrightarrow M$ is the inclusion. In
general, the symplectic quotient is a symplectic stratified space
\cite{ls}. Suppose that the original manifold $(M, \Omega)$ is
K\"ahler, and that the action preserves both the K\"ahler metric
and the symplectic form, hence, the complex structure. Then the
symplectic quotient $M{\sslash}G$ is K\"ahler with K\"ahler form
$\ow$ \cite{gs2, hitchin}.

Products of compact non-Lefschetz symplectic manifolds with a
two-sphere rotated about a fixed axis are trivial examples such that
the existence of a non-Lefschetz symplectic quotient implies the original
manifold is not Lefschetz.
However, Lin \cite{lin} constructed examples of
Hamiltonian $S^1$-manifolds $(M, \Omega)$ that are strong Lefschetz but their
symplectic quotients $(M{\sslash}S^1, \ow)$ are not. In this paper, we
generalize Lin's result.

Let $\ow_{\C^n}=\frac{i}{2}\sum dz_j\wedge d\bar{z}_j=\sum dx_j\wedge
dy_j$ be the standard symplectic form on \(\C^n.\)
Let \(G\) be a compact Lie group with an \(n\)-dimensional
unitary representation $\tau\colon G \to \U(n)$.
Then \((\C^n,
\omega_{\C^n})\) is a Hamiltonian \(G\)-manifold.

\begin{theorem} \labell{T:theorem}
Suppose \(\C^n{\sslash}G\) consists of a single point.
Then there exists a compact connected Hamiltonian \(G\)-manifold $(M, \Omega)$
that is strong Lefschetz,
but its symplectic quotient $(M{\sslash}G, \ow)$ is not.
\end{theorem}

As an application, let \(m, k = 2, 3, 4, \dots\) and \(G\) be either
\(\Su(m)\) or \(\Sp(k)\).  Then the moment map
for the standard \(G\)-representation on $\C^n$, for $n=m$ and $n=2k$,
respectively, is
\[
\psi\colon \C^n \to \fg \cong \fg^* \qquad \text{given by}\qquad
\psi(z) = \frac{i}{2}zz^*,
\]
where we identify the Lie algebra \(\fg\) with its dual $\fg^*$
via the inner product $\la A, B\ra= \text{trace}(A^*B)$.  A direct
calculation shows that \(\psi^{-1}(0) = \{0\}\).
Theorem~\ref{T:theorem} then implies

\begin{theorem}\labell{T:app}
Let \(n > 1\) and let
\(G\) be either \(\Su(n)\) or \(\Sp(n)\).  Then
there exists a compact connected Hamiltonian \(G\)-manifold $(M, \Omega)$
that is strong Lefschetz, but its symplectic quotient $(M{\sslash} G, \ow)$ is
not.
\end{theorem}


In Section \ref{S:construction}, we spell out the construction for
the manifold $(M, \Omega)$ in Theorem \ref{T:theorem}. In
Section \ref{S:cohomology}, we compute the Lefschetz maps.
In our proof, we use a result of Gompf \cite{gompf} as our symplectic
quotient:

\begin{theorem} [Gompf] \labell{L:gompf}
There exists a
4-dimensional symplectic manifold \((B,\omega)\) such that the
Lefschetz map
$\ow \colon \H^1(B) \to \H^3(B)$
is the zero-map. Moreover, there exists an integral class
\(c \in \H^2(B)\) such that
$c \colon \H^1(B) \to \H^3(B)$
is an isomorphism.
\end{theorem}


The strong Lefschetz property also has an equivalent
description in terms of symplectic harmonic forms.  For details, see
\cite{jeanluc, mathieu, yan}.
 The Hamiltonian
\(G\)-manifolds \((M,\Omega)\) we construct here have the property
that each de Rham cohomology class contains a symplectic harmonic
representative while the same is not true for \((M{\sslash}G,
\ow)\).

\section{Construction}\labell{S:construction}


Let $G$ be a compact Lie group acting on $(\C^n, \ow_{\C^n})$ by a
unitary representation with the property that the symplectic
quotient $\C^n{\sslash}G$ is a point. Extending the action to
$\C^n \times \C$ trivially on the last factor, it induces a
Hamiltonian $G$-action on the complex projective space $\PP(\C^n
\times \C)=\CP^n$. Let $F = \CP^n$ and $\omega_F$ the Fubini-Study
symplectic form on F scaled by $\epsilon > 0$. The symplectic
quotient $F{\sslash}G$ under the induced $G$-action is again a
point.

\begin{remark}\labell{R:cut}
This projective space $(F, \ow_F)$ can be constructed by
symplectic reduction. Let $S^1$ denote the circle group and we
identify both the Lie algebra of $S^1$ and its dual with $\R$. The
multiplication of $S^1$ on $\C^n \times \C$ is a Hamiltonian
action with moment map $\phi\colon \C^n\times \C \to \R$ given by
$$\phi(z,w)=\frac{1}{2}(-\|z\|^2-|w|^2+ \text{constant}).$$ If we
choose the constant to be $\epsilon$, the symplectic quotient is
$(F, \ow_F)$. This is the same as taking a closed ball in $\C^n$
of radius $\epsilon^{1/2}$ centered at zero and reduce the
boundary since
\[
\begin{split}
\phi^{-1}(0)/S^1 &= \left\{\|z\|^2+|w|^2=\epsilon\right\}/S^1 \\
&\cong \left\{\|z\|^2 =\epsilon\right\}/S^1 \sqcup \left\{\|z\|^2
< \epsilon\right\}.
\end{split}
\]
\end{remark}

Let $(B, \ow)$ be a symplectic manifold and let
\(\pi\colon P \to B\)
be a principal $S^1$-bundle over \(B\) with a given Chern class
\(c\) and \(\theta\) a connection form on \(P\). With the Hamiltonian
$S^1$-action
on $(F, \ow_F)$ induced from the multiplication on $(\C^n, \ow_{\C^n})$, we
form the associated bundle $P \times_{S^1} F$.

\begin{lemma}\labell{L:coupling}
For $\epsilon$ sufficiently small, the associated bundle $P
\times_{S^1} F$ is a symplectic manifold with fiber $(F, \ow_F)$.
\end{lemma}



\begin{proof}

This is done by the technique of minimal coupling
\cite{sternberg, weinstein2}.

First identify both the Lie algebra of $S^1$ and its dual with
$\R$ and consider the space \(P \times \R\). The minimal coupling
form
\[
\pi^* \omega + d \la \text{pr}_2,\theta \ra
\] is nondegenerate on $P\times \{0\}$, and therefore is symplectic
on a small \(\delta\)-neighborhood \(P \times (- \delta,
\delta)\). Let $I = (-\delta, \delta)$.
The \(S^1\)-action on \(P \times I\) given by
$a \cdot (p,\eta) = (pa^{-1}, \text{Ad}^*(a)\eta)$
is Hamiltonian
with minus the projection onto
the second factor
\[
-\text{pr}_2\colon P \times I \to \R
\]
being its moment map. The symplectic quotient
$(P\times I) {\sslash} S^1$ is $(B,\ow)$.

Second consider the product space $P \times I \times F$. Let
$\Phi_F\colon F \to \R$ be the $S^1$-moment map on $F$. The
diagonal $S^1$-action on $(P\times I) \times F$ is Hamiltonian
with moment map given by
\[
\Phi(p,\eta,f) = \Phi_F(f) - \eta.
\]
If $\epsilon < \delta$, the image of $F$ under $\Phi_F$ is contained
in $I$, and $\Phi^{-1}(0) \cong P \times
F$. So the symplectic quotient is
\[
\Phi^{-1}(0)/S^1 = P \times_{S^1} F
\]
with the symplectic form
\[
\Omega([p,f]) = \pi^* \omega(p) + \la \Phi_F(f),\Theta(p) \ra +
\ow_F(f),
\]
where $\Theta = d\theta - \frac{1}{2}[\theta, \theta]$ is the
curvature form of \(\theta\). Since \(\theta\) is a connection
1-form on the principal $S^1$-bundle \(P\) with Chern class \(c\),
it follows that $\pi^*c = -[{\Theta}]$.

\end{proof}

Finally, since the $G$-action on $\C^n$ we start out with commutes
with the $S^1$-multiplication, the induced Hamiltonian actions on
$F$ also commute. The $G$-moment map on $F$ is $S^1$-invariant.
Hence the manifold \((P \times_{S^1} F, \Omega)\) inherits a
fiberwise \(G\)-action with moment map
\[\Psi\colon P
\times_{S^1} F \to \fg^*\] induced by the $G$-moment map on the
fiber. Reduction can be carried out fiber by fiber. Since
$F{\sslash}G$ is a point, by construction, $\Psi^{-1}(0)/G = (B,
\ow)$.

\begin{lemma}
Let $M=P\times_{S^1} F$ and $\Omega$ be as in the proof of Lemma
\ref{L:coupling}. Then $(M, \Omega)$ is a Hamiltonian $G$-manifold
with $M{\sslash}G = (B, \ow)$.
\end{lemma}


For our purpose,
we use the symplectic manifold $(B, \ow)$ and the Chern class $c$
prescribed in Theorem~\ref{L:gompf}.


\section{Strong Lefschetz property}\labell{S:cohomology}

This section proves that the Hamiltonian \(G\)-manifold
\((M, \Omega)\) constructed in the previous section
is strong Lefschetz. Namely, we show that
the map
\[
\Omega^{k}\colon \H^{n+2-k}(M) \to \H^{n+2+k}(M)
\]
is an isomorphism for every integer \(0 \le k \le n+2\). By
Poincar\'e Duality, it suffices to show that the map $\Omega^k$ is
injective. This holds for $k=n+2$ by the nondegeneracy of the
symplectic form and also for $k=0$. It remains to show injectivity
for $1 \leq k \leq n+1$.

A typical fibre of $(M, \Omega)$ is $(F, \ow_F)$, which is
isomorphic to \(\CP^n\) with a scaled Fubini-Study
symplectic form.  To simplify notations, let
\[
x = \pi^* \omega, \ \ y = \la \Phi_F, \Theta
\ra, \ \ z = \ow_F, \ \ u = \la \Phi_F, \Theta
\ra + \ow_F = y + z.
\]
When there is no confusion, we use the same symbol to denote a
closed form and its cohomology class.  For example, \(x\) and
\(u\) are both closed and we use the same symbol for both the form
and its cohomology class in the subsequent computation.

In this notation,
\(\H^*(F)\) is generated by the class
\(u|_F \in \H^2(F)\).  Hence \(\H^*(M)\) as a
\(\H^*(B)\)-module is freely generated by \(\{1, u, u^2, \dots,
u^n\}\) by the Leray-Hirsch Theorem (see, for example, \S 5 of \cite{bt}):

\begin{theorem}[Leray-Hirsch]\labell{LH}
Let $M$ be a fiber bundle over $B$ with fiber $F$.
The cohomology \(\H^*(M)\) is a free \(\H^*(B)\)-module
generated by \(\{e_1,\cdots,e_r\} \subset \H^*(M)\) if the restriction
of \(\{e_1,\cdots,e_r\}\) to \(F\) generates \(\H^*(F)\).
\end{theorem}


For the rest of the paper, we identify the elements in
$\H^*(B)$ and $\H^*(F)$ with their images in $\H^*(M) \cong \H^*(B)\otimes
\H^*(F)$. By Theorem \ref{LH}, a class $\alpha \in \H^k(M)$
can be written as
\begin{align*}
\alpha  =
\begin{cases}
b_4u^{(k-4)/2}+b_2u^{(k-2)/2}+b_0u^{k/2}, & \text{if $k$ is even,}\\
b_3 u^{(k-3)/2}+b_1u^{(k-1)/2}, & \text{if $k$ is odd,}
\end{cases}
\end{align*}
where $b_j$ denotes both a class in $\H^j(B)$ and its lifting in
$\H^j(M)$.

In particular,
\[u^{n+1} =\beta_4 u^{n-1} + \beta_2 u^n
= \frac{n(n+1)}{2}y^2z^{n-1}+(n+1)yz^n \] for some \(\beta_2 \in
\H^2(B)\) and \(\beta_4 \in \H^4(B)\).

Let \(h \in \H^2(B)\).  Then
\[
\begin{split}
u^{n+1} h &= (\beta_4 u^{n-1}+\beta_2 u^n)h = \beta_2 h u^n =
\beta_2 h z^n \\
&= \left(\frac{n(n+1)}{2}y^2z^{n-1}+(n+1)yz^n \right) h = (n+1)yhz^n.
\end{split}
\]
The form \(z^n\) restricted to each fibre of the bundle projection \(M
\to B\) is the volume form on
\((F,\omega_F)\). Integrating
along the fibre, we obtain
\[
\int_M \beta_2 h z^n = V(\omega_F) \int_B \beta_2 h,
\]
where \(V(\omega_F) = \frac{\pi^n \epsilon^{n}}{\Gamma(n+1)}\) is
the volume and $\Gamma$ is the Euler $\Gamma$-function. On the
other hand,
\[
\int_M (n+1)yh z^n = \int_M (n+1) \la \Phi_F, \Theta \ra h z^n =
\int_M -(n+1) \la \Phi_F, \pi^*c \ra h z^n.
\]
In this case, along each fibre, the volume form \(\omega_F^n\) is
scaled by the moment map \(\Phi_F\). The volume of \(S^{2n-1}\)
with radius \(r\) is \(\frac{2\pi^nr^{2n-1}}{\Gamma(n)}\). By
Remark \ref{R:cut}, it follows that
\[
\begin{split}
\int_M (n+1)yh z^n&=\left[\int_0^{\epsilon^{1/2}} \frac{r^2}{2}\cdot
\frac{2\pi^n
r^{2n-1}}{\Gamma(n)}dr\right]  \int_B (n+1)ch \\
&= \frac{\pi^n\epsilon^{n+1}}{(2n+2)\,\Gamma(n)}\int_B (n+1)ch\\
&=\frac{\pi^n \epsilon^{n+1}}{2\, \Gamma(n)}\int_B ch\,.
\end{split}
\]
Hence $\beta_2 = \frac{1}{2}n \epsilon c$. Note
that \(\beta_2\) depends on \(\epsilon\) which can be chosen
anywhere between \(0\) and \(\delta\).

Similarly,
\[
u^{n+2}=(\beta_4+\beta_2^2)u^n=\frac{(n+1)(n+2)}{2}y^2z^n.
\]
Integrating, we get
\[
\begin{split}
\int_M u^{n+2} &= \int_M (\beta_4+\beta_2^2)u^n = \int_M
(\beta_4+\beta_2^2)z^n
= V(\ow_F)\int_B (\beta_4+\beta_2^2)\\
&= \int_M \frac{(n+1)(n+2)}{2}y^2z^n =
\frac{(n+1)\pi^n\epsilon^{n+2}}{8\, \Gamma(n)}\int_B c^2.
\end{split}
\]
Hence $\beta_4 = \frac{1}{8}n(1-n)\epsilon^2 c^2$.

Now we are ready to compute the kernel of the Lefschetz map
\[\Omega^k\colon \H^{n+2-k}(M) \to \H^{n+2+k}(M).\] Let \(\alpha \in
\H^{n+2-k}(M)\). We will show that $\Omega^k\alpha =0$ implies
that $\alpha=0$.

With the notations explained earlier, we have
\begin{align*}
\Omega^k &= (x+u)^k \\
&=\begin{cases}
x+u, & \mbox{ if } k=1, \\
\frac{k(k-1)}{2}x^2 u^{k-2} + kx u^{k-1}
+ u^k, & \mbox{ if } 2\leq k \leq n,\\
\left(\frac{n(n-1)}{2}x^2 + \beta_4\right)u^{n-1} + ((n+1)x +
\beta_2)u^n, & \mbox{ if } k=n+1.
\end{cases}
\end{align*}

If \(k=n+1\), we have $\alpha = b_1$. Hence,
\[
\Omega^{n+1} \alpha = ((n+1) x + \beta_2)b_1 u^{n}.
\]
If $\alpha$ is in the kernel, that is, $\Omega^{n+1} \alpha = 0$,
then
\[
((n+1) x + \beta_2)b_1
= 0.
\]
Since $(n+1) x + \beta_2 = (n+1) x + \frac{1}{2} n \epsilon c$
defines an injective map from $\H^1(B)$ to $\H^3(B)$ for a generic
\(\epsilon\), we conclude that $b_1=0$. And therefore
$\ker(\Omega^{n+1})= \{0\}$ for a generic $\epsilon$.

If \(k=n\), we have $\alpha = b_2+ b_0 u$.
Hence,
\[
\begin{split}
\Omega^n \alpha &=
\left(b_2 nx+ b_0 \frac{n(n-1)}{2} x^2 \right) u^{n-1}
+ (b_2+ b_0 n x) u^n
+ b_0 u^{n+1}\\
&= \left(b_2 nx + b_0\frac{n(n-1)}{2}  x^2+ b_0\beta_4 \right) u^{n-1}
+ (b_2 + b_0 nx+ b_0\beta_2) u^n.
\end{split}
\]
If $\alpha$ is in the kernel, that is, $\Omega^{n} \alpha = 0$,
then
\[
\left\{
\begin{array}{lll}
b_2 n x + b_0 (\frac{n(n-1)}{2}x^2+\beta_4) & = & 0,\\[2ex]
b_2 +b_0 (nx+\beta_2) & = & 0.
\end{array}
\right.
\]
Assume that $b_0 \neq 0$. Since
$b_2 = -b_0(nx+ \beta_2)$, we get
\[
\begin{split}
0 &= -b_0(nx+ \beta_2) n x + b_0
\left(\frac{n(n-1)}{2}x^2+\beta_4\right)\\
& = \left(\beta_4 - n\beta_2 x - \frac{n(n+1)}{2}x^2 \right)b_0.
\end{split}
\]
This implies that
\[
\frac{n(1-n)\epsilon^2}{8} c^2 - \frac{n^2\epsilon}{2} cx -
\frac{n(n+1)}{2}x^2  = 0.
\]
However the left hand side is lifted from an element in $\H^4(B)$
and can be made non-zero with a generic $\epsilon$. This is a
contradiction. Hence for a generic \(\epsilon\), we must have
$b_0=0$ and $b_2=-b_0(\beta_2 + nx)=0$. We conclude that
\(\ker(\Omega^n) = \{0\}\) for a generic \(\epsilon\).

If \(1 \leq k< n \) and $n-k$ is even, we have
\[
\alpha = b_4
u^{(n-k-2)/2}+  b_2 u^{(n-k)/2} + b_0 u^{(n-k+2)/2}.
\]
Hence,
\[
\begin{split}
\Omega^k \alpha &=
\left(b_4 + b_2 k x + b_0 \frac{k(k-1)}{2}x^2\right) u^{(n+k-2)/2}\\
&\qquad+ (b_2 +b_0 k x)  u^{(n+k)/2}
+ b_0 u^{(n+k+2)/2}.
\end{split}
\]
Suppose \(\Omega^k \alpha = 0\).  Then
\[
\left\{
\begin{array}{lll}
b_4 + b_2 k x + b_0 \frac{k(k-1)}{2}x^2 & = & 0,\\[2ex]
b_2 + b_0 k x & = & 0,\\[2ex]
b_0 & = & 0.
\end{array}
\right.
\]
This implies that $b_0 = b_2 = b_4 = 0$.
Hence \(\ker(\Omega^k) = \{0\}\).

If \(1 \leq k< n \) and $n-k$ is odd, we have
\[
\alpha = b_3
u^{(n-k-1)/2}+  b_1 u^{(n-k+1)/2}.
\]
Hence,
\[
\Omega^k \alpha = (b_3 +b_1 k x)
u^{(n+k-1)/2}+  b_1 u^{(n+k+1)/2}.
\]
Suppose \(\Omega^k \alpha = 0\).  Then
\[
\left\{
\begin{array}{lll}
b_3 + b_1 k x & = & 0,\\[2ex]
b_1 & = & 0.
\end{array}
\right.
\]
This implies that $b_1 = b_3 = 0$.
Hence \(\ker(\Omega^k) = \{0\}\).

Hence we conclude that for a generic \(\epsilon\), the Lefschetz map
\(\Omega^k\) is an isomorphism for all \[1 \leq k \leq n+1.\]
This completes the proof
of Theorem \ref{T:theorem}. Our manifold $(M, \Omega)$ is strong
Lefschetz while its symplectic quotient $(B, \ow)$ is not.

\end{document}